\documentclass{ifacconf}

\usepackage{amsmath}
\usepackage{graphicx}      
\usepackage{natbib}        
\begin{document}
\begin{frontmatter}

\title{Identification and Online Updating of Dynamic Models for Demand Response of an Industrial Air Separation Unit}


\author[First]{Calvin Tsay}
\author[Second]{Yanan Cao}
\author[Second]{Yajun Wang}
\author[Second]{Jesus Flores-Cerrillo}
\author[First]{Michael Baldea}

\address[First]{McKetta Department of Chemical Engineering, The University of Texas at Austin, Austin, TX 78712 USA (e-mail: calvint@utexas.edu, mbaldea@che.utexas.edu).}
\address[Second]{Smart Operations, Center of Excellence (COE), Linde, Tonawanda, NY 14150 USA (e-mail: yanan.cao@linde.com, yajun.wang@linde.com, jesus.flores-cerrillo@linde.com).}

\begin{abstract}                
Demand-response operation of air separation units requires frequent changes in production rate(s), and scheduling calculations must explicitly consider process dynamics to ensure feasibility of the solutions. 
To this end, scale-bridging models (SBMs) approximate the scheduling-relevant dynamics of a process and its controller in a low-order representation.  
In contrast to previous works that have employed nonlinear SBMs, this paper proposes linear SBMs, developed using time-series analysis, to facilitate online scheduling computations. 
Using a year-long industrial dataset, we find that compact linear SBMs are suitable approximations over typical scheduling horizons, but that their accuracies are unpredictable over time. 
We introduce a strategy for online updating of the SBMs, based on Kalman filtering schemes for online parameter estimation. The approach greatly improves the accuracy of SBM predictions and will enable the use of linear SBM-based demand-response scheduling in the future.  

\end{abstract}


\end{frontmatter}

\section{Introduction}
Cryogenic air separation units (ASUs) are strong candidates for demand-response operation. 
The primary cost of ASU operation is associated with electricity usage, as compressors driven by large electric motors are used to bring feed air to the requisite operating pressures.
As a result, process economics can be greatly improved by proactively scheduling production in response to (predicted) fluctuations in electricity price---which may become more extreme with increasing adoption of renewable-based generation.
Moreover, the products of an ASU comprise purified components of air, which can be liquefied and stored.
The above factors motivate the development of production schedules focused on ``load shifting'': increasing production when electricity prices are low and storing excess products in liquefied form, which can later be used to satisfy product demand(s) when electricity prices are high.

There are several modeling considerations for demand-response scheduling of ASUs. 
To fully exploit time-varying electricity prices, which change hourly (or more frequently), production schedules must consider hourly (or shorter) scheduling intervals.
Using such scheduling intervals requires that ASU dynamics be considered explicitly to ensure that production schedules are both dynamically feasible and economically optimal. 
Embedding a representation of process dynamics in production scheduling models requires a multi-scale approach that bridges the (faster) time scales of process dynamics/control with those of production scheduling \citep{jamaludin2017approximation}. 
However, such models must be used to optimize process operation over relatively long time horizons, yet still accurately represent complex process behavior over all relevant time scales. 
Several reduced-order modeling approaches \citep{cao2016dynamic, caspari2020wave, schafer2019reduced} have been proposed for ASUs, mostly derived using physical insight and engineering expertise. 
Data-driven, system identification-based methods have also been described \citep{pattison2016optimal, tsay2019110th}.  

In particular, our recent works \citep{tsay2019optimal,  tsay2020integrating} demonstrated that \textit{nonlinear} scale-bridging models (SBMs), or low-order representations of closed-loop process dynamics, can be identified from historical ASU operating data. 
Production scheduling is then cast as a dynamic optimization problem using SBMs to represent the dynamics of the ASU and its control system. 
Nevertheless, developing nonlinear SBMs requires significant system identification effort.
Furthermore, SBM approximations may not remain accurate over extended periods of time, and the parameters may periodically require re-estimation.
This degradation of accuracy can be attributed to the SBM not capturing all process/controller dynamics, to system identification only being performed on certain regions of the input space, and/or to changes in the process dynamics and control system over time.

Motivated by these challenges, in this paper we introduce an approach for online updating of SBMs used for production scheduling. 
Using data recorded over a year of routine operation of an industrial cryogenic ASU, we show that the proposed Kalman filter-based strategy enables the use of low-order, \textit{linear} SBMs (in the form of simple time-series models) over typical scheduling time horizons. 
We show that the online update strategy maintains model accuracy over the relevant time scales, which would otherwise degrade quickly over time.

\section{Process and Data Description}

In a cryogenic ASU, the feed stream passes through a feed air compressor (FAC) before being sent to the heat exchanger and separation columns. 
Further details regarding ASU processes and their control systems can be found in, e.g., \cite{caspari2020integration, tsay2019optimal}.
The FAC accounts for the majority of the electricity consumed by the ASU, and demand-response scheduling calculations therefore rely on accurate predictions of the FAC power consumption, denoted as $\dot{W}_\mathrm{FAC}$, as a function of production rate and other factors. 
Given its importance, in this work we focus on the modeling of $\dot{W}_\mathrm{FAC}$. 
Note that the techniques described here can be applied to all scheduling-relevant variables needed to formulate a complete production scheduling optimization problem \citep{pattison2016optimal, tsay2020integrating}.

\begin{figure}
	\begin{center}
		\includegraphics[width=\columnwidth]{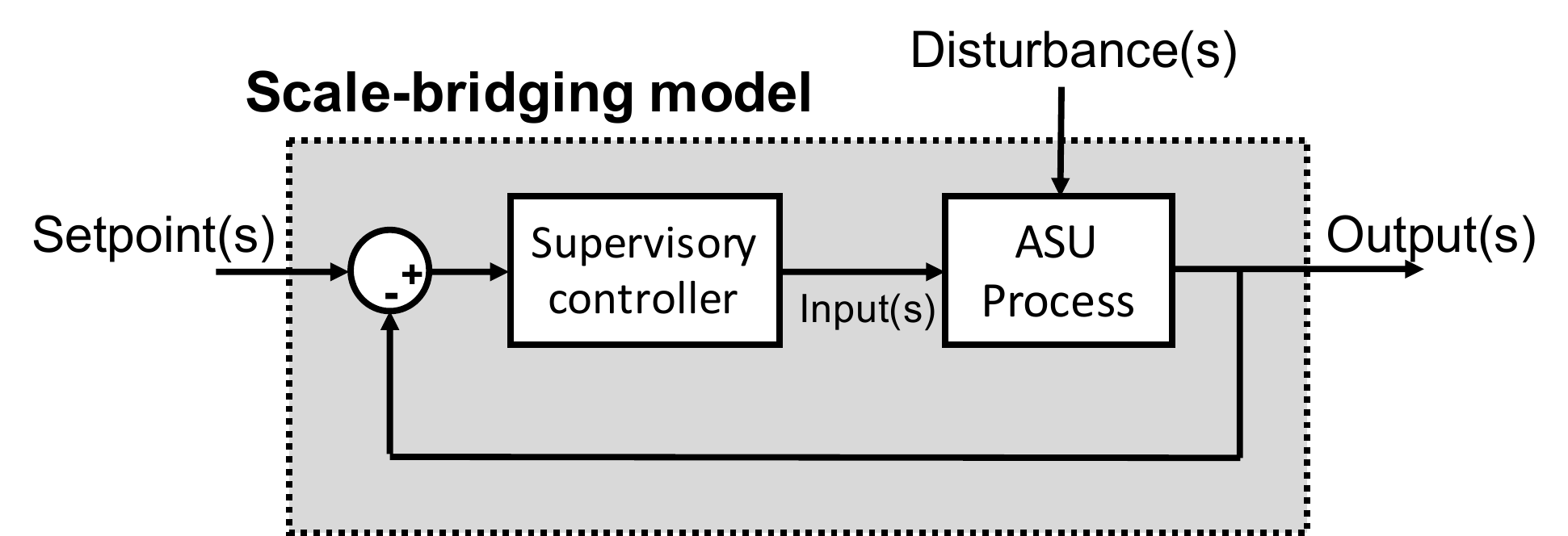}    
		\caption{Scale-bridging model (SBM) concept. The shaded area reflects the system whose input-output dynamics are represented by the SBM.}
		\label{fig:sbm}
	\end{center}
\end{figure}

Introduced by \cite{du2015time}, SBMs aim to represent the closed-loop, input-output response of process variables, with the inputs being controller setpoints and/or measured disturbance variables (Figure \ref{fig:sbm}). 
As a result, system identification is performed on the \textit{closed-loop} response of the process, in contrast to  conventional \textit{open-loop} system identification (where the input data are process inputs). 
The ASU operates under multivariable model predictive control (MPC), with eight operator setpoints available. 
Ambient temperature ($T$) acts as a measured disturbance variable. 
Data for setpoints, $T$, and $\dot{W}_\mathrm{FAC}$ were recorded at one-minute intervals over one year of routine operation. 

An unexpected feature of the dataset was that data for six of the setpoints are strongly correlated. 
Principal component analysis (PCA) for these six setpoints reveals that 90\% of the variation is explained by the first two principal components. 
Due to this correlation, the individual effects of these six setpoints cannot be isolated using the available data, and they are replaced with the first two principal components ($\phi_1$, $\phi_2$) to reduce model dimensionality. 
Of the remaining two setpoints, one is uncorrelated ($<$1\%) with $\dot{W}_\mathrm{FAC}$ and is not treated as an SBM input. The other, denoted as $\mathrm{SP}_1$, is treated as an input. In total, four inputs were used to construct the SBM: $\phi_1$, $\phi_2$, $\mathrm{SP}_1$, and $T$.

All input and output variables are normalized and filtered to reduce the amount of high-frequency noise, and to protect the confidential industrial data. 
SBM input variables are filtered using a low-pass filter with a cutoff of 1 h$^{-1}$, while the output variable ($\dot{W}_\mathrm{FAC}$) is treated using a Savitzky-Golay filter with a window of 15 samples.

\section{Scale-Bridging Model Development}
\label{sec:sbm}
\subsection{Time-Series Modeling}

\begin{figure}
	\begin{center}
		\includegraphics[width=\columnwidth]{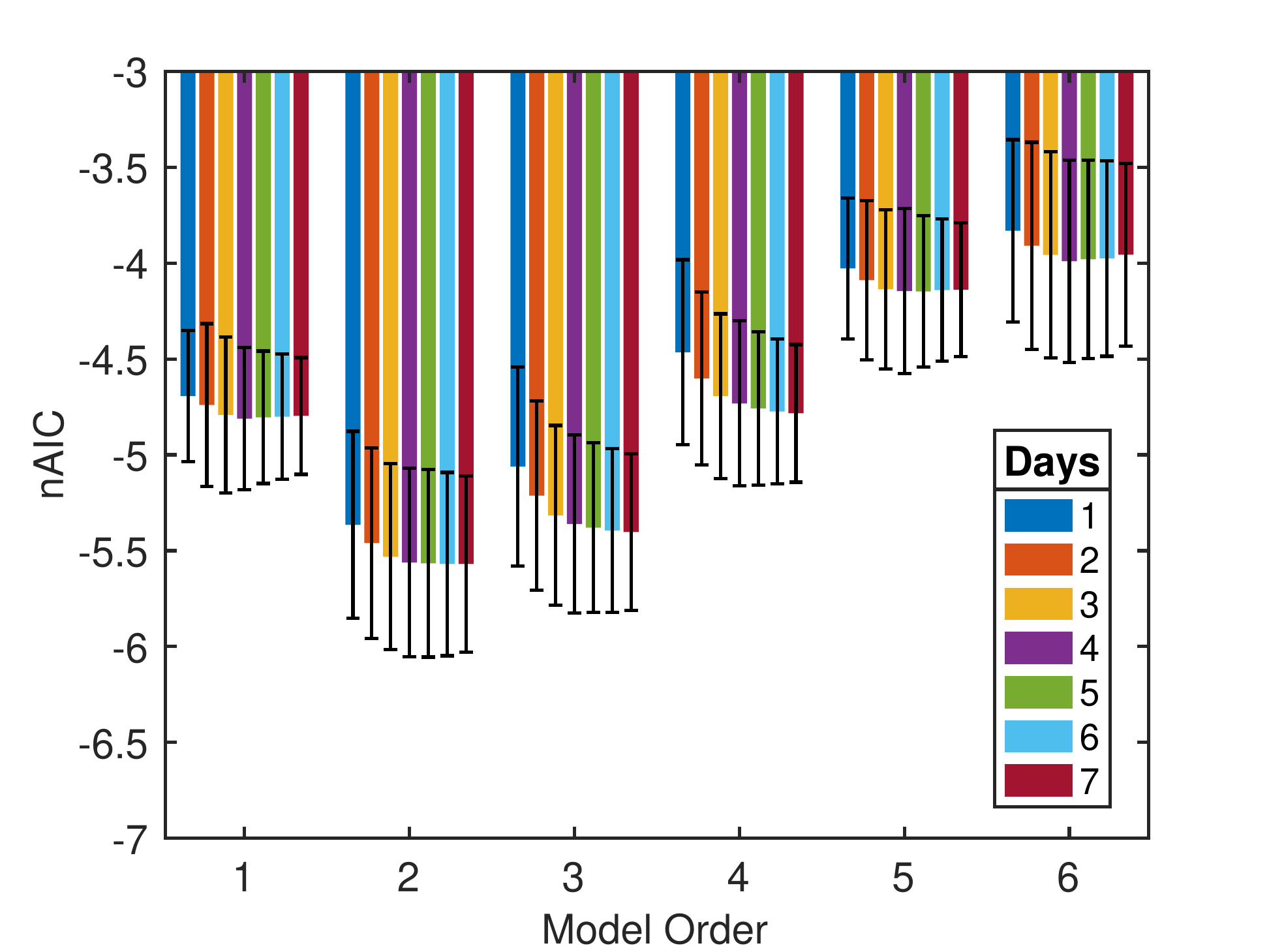}    
		\caption{Normalized AIC for varying model order $N$ and number of days in $D_\mathrm{train}$. Error bars denote one standard deviation.}
		\label{fig:aic}
	\end{center}
\end{figure}

Our previous work \citep{tsay2019optimal} showed that \textit{nonlinear} SBMs (i.e., Hammerstein-Wiener models) can accurately represent process dynamics over the span of several months. 
While this avoids frequent re-fitting of the SBM parameters, demand-response scheduling horizons are typically only several days in length, suggesting that even simpler model forms can be used. 
Therefore, in this work, we study \textit{linear} SBMs, which simplify both the system identification step and potentially the later dynamic-optimization-based scheduling step. 
For example, \cite{kelley2018milp} formulated the latter as a MILP by linearizing SBM dynamics, while integer variables were introduced to handle SBM nonlinearities.
A linear SBM could eliminate the need for integer variables, making the scheduling problem an LP.
We note that \cite{dias2018simulation} proposed embedding linear, \textit{open-loop} dynamic models along with the process MPC in scheduling calculations.

\begin{figure*}
	\begin{center}
		\includegraphics[width=\textwidth, trim=3cm 0cm 3cm 0.8cm, clip]{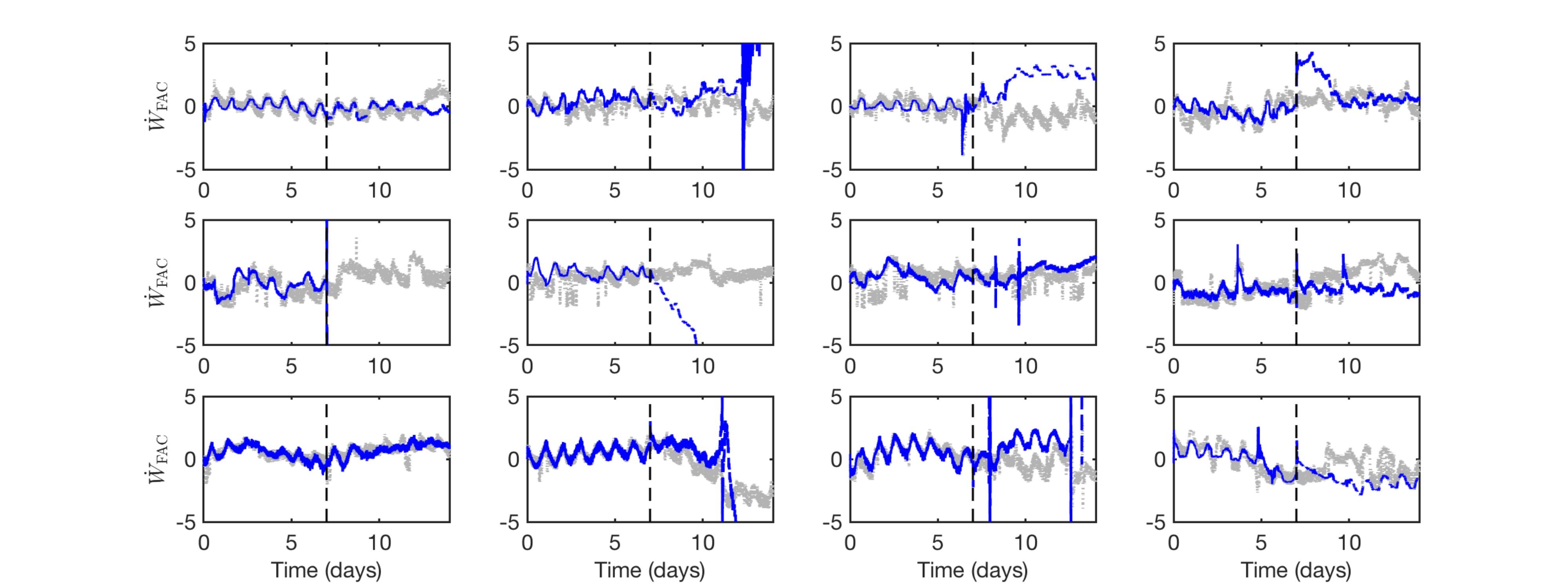}    
		\caption{SBM predictions for each month. The left-hand side of each plot (left of the dashed line) shows  predictions for training data, while the right-hand side shows predictions over the next seven days. Recorded data are shown in gray.}
		\label{fig:sbm_pred}
	\end{center}
\end{figure*}

We consider linear SBMs in the form of an autoregressive with extra inputs (ARX) time-series model:
\begin{equation}
y_t + \sum_{j=1}^{N} a_j y_{j-i} = \sum_{i=1}^{M} \sum_{j=1}^{N} b_{i,j} u_{i, t-j} + e_t \label{eq:arx}
\end{equation} 
where $y_t$ is the output value at time $t$, and $u_{i,t}$ is the value of the $i^\mathrm{th}$ input ($i = 1,...,M$) at time $t$. 
The coefficients $a_j$ are autoregressive parameters, $b_{i,j}$ are the regressed coefficients for the $i^\mathrm{th}$ input, and $e$ is the model error. 
Note that the linear time-series model \eqref{eq:arx} can easily be converted to a continuous-time form, such as a transfer function or state-space model. 
In practice, the number of autoregressive terms (number of poles) and the number of coefficients for each input (number of zeros) are independent; however, for simplicity we treat both of these as a single model order, denoted as $N$.

We consider each month of the one-year dataset separately, yielding 12 datasets in total. 
SBM models of the form  \eqref{eq:arx} are trained using the first $D_\mathrm{train}$ days of data for each month. 
Figure \ref{fig:aic} shows the normalized Akaike Information Criterion (nAIC) for several different model orders $N$ and values of $D_\mathrm{train}$. 
Models with $N$=2 and $N$=3 have similar nAIC, particularly as $D_\mathrm{train}$ increases. 
Thus, the rest of this work considers SBMs of the form \eqref{eq:arx} with $N$=3, which provide a good tradeoff between model size and accuracy.

\subsection{System Identification Results}

\begin{figure}
	\begin{center}
		\includegraphics[width=0.8\columnwidth, trim=1cm 2cm 1cm 1.5cm, clip]{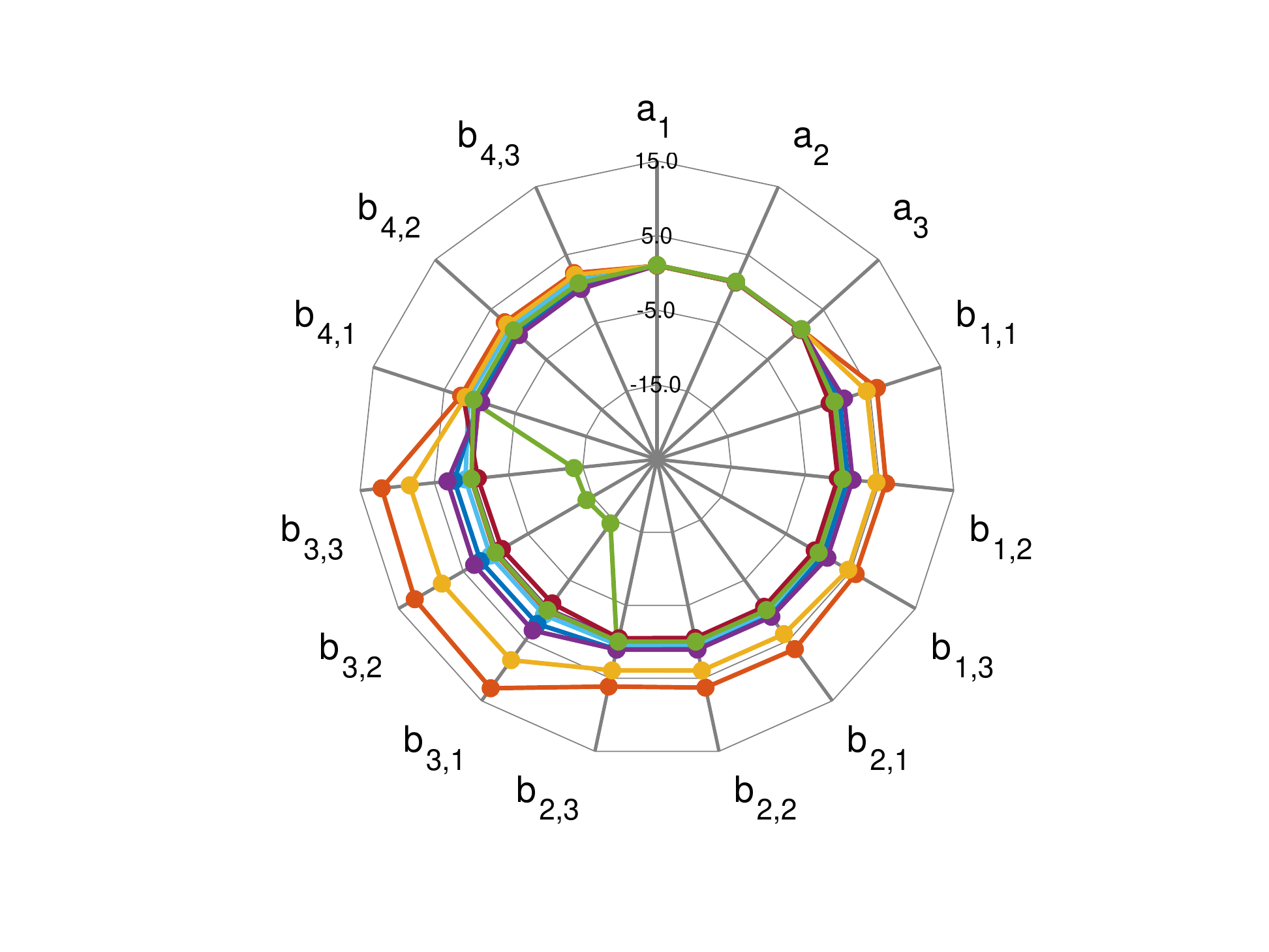}    
		\caption{Scaled SBM parameter values. Each plot/ring represents one SBM (one for each of the 12 months).}
		\label{fig:params_radar}
	\end{center}
\end{figure}

Figure \ref{fig:sbm_pred} shows the predictions of the trained SBMs on seven days of training data, as well as the next seven days. The mean squared error (MSE) for predictions on training data (left side of each plot) for the twelve datasets is $0.47\pm0.27$ (mean $\pm$ standard deviation), while the MSE for predictions on new data (right side of each plot) is $51.67\pm94.78$. 
The simple, third-order ARX model can accurately represent the closed-loop process dynamics over one week of training data, but its accuracy is unpredictable over time; some SBMs in Figure \ref{fig:sbm_pred} make good predictions on new data, while others quickly lose accuracy. 

The fitted parameter values for the trained SBMs are shown in Figure \ref{fig:params_radar}. The plotted values are scaled by the respective mean over the 12 SBMs. 
While the values of some parameters (e.g., $b_{3,j}; j=1,2,3$) varied significantly among SBM models, the values of $a_j$ and $b_{4,j}$ were nearly constant. 
Interestingly, this data-driven approach reveals that the system poles (i.e., closed-loop time constants) do not change significantly over the 12 monthly datasets. 
The dependence of $\dot{W}_\mathrm{FAC}$ on $u_4$ (ambient temperature $T$) also remains relatively stable. 
These observations confirm the physical intuition that the response of the process to the controller setpoints $\{u_1, u_2, u_3 \}  = \{ \phi_1, \phi_2, \mathrm{SP}_1 \}$ is nonlinear and/or time-varying (e.g., as the control system is manipulated by process operators), and cannot be represented accurately using a linear dynamic model with fixed parameters. 
On the other hand, the poles and the response to changes in ambient temperature are characteristic to the process itself and do not change over the time span of the dataset.

\begin{figure}
	\begin{center}
		\includegraphics[width=\columnwidth]{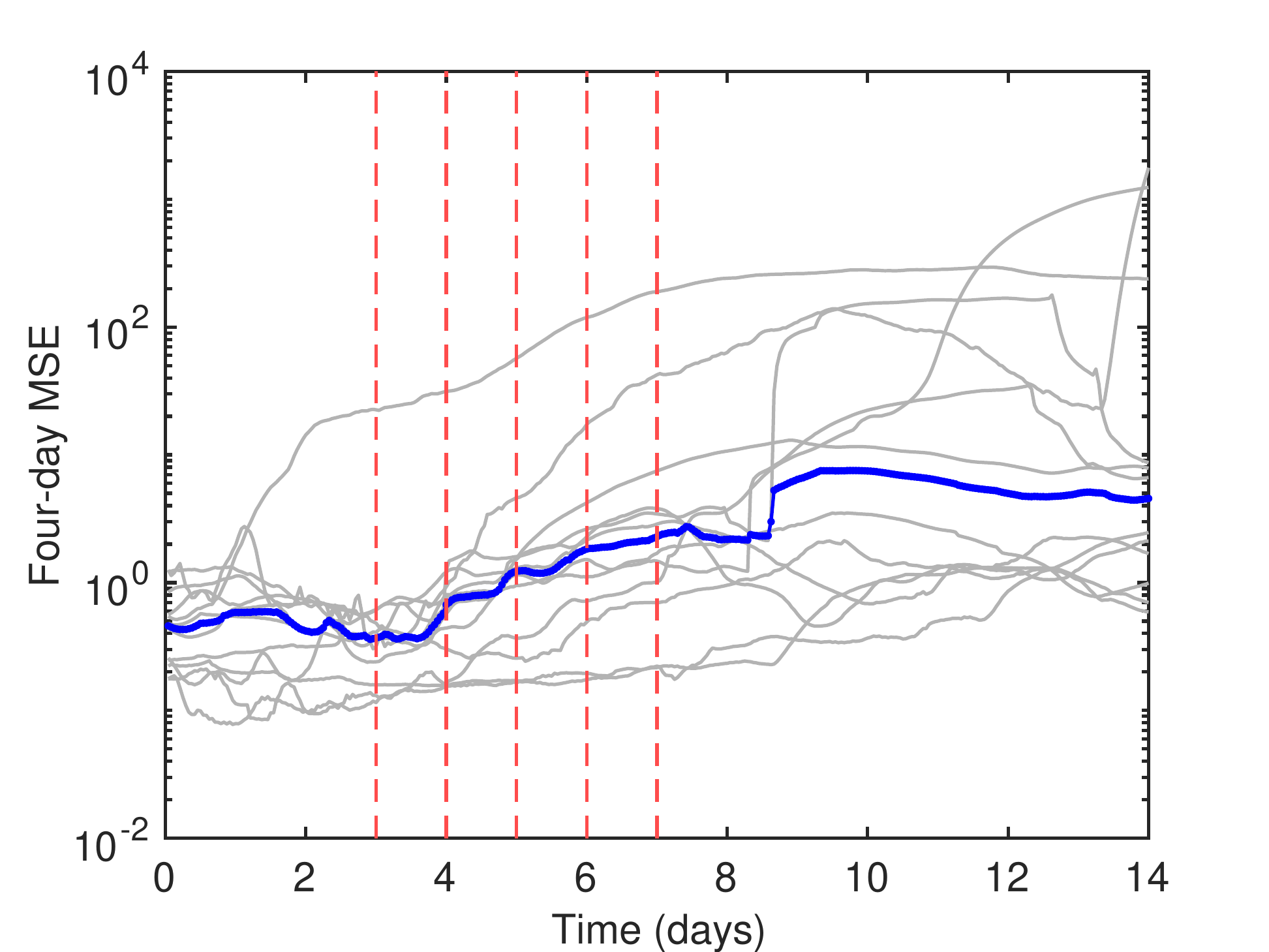}    
		\caption{SBM accuracy over a four-day moving horizon. Each gray line depicts one SBM (from one monthly dataset), while the blue line indicates the median values. The dashed vertical lines mark days 3--7.}
		\label{fig:sbm_mse}
	\end{center}
\end{figure}

While Figure \ref{fig:sbm_pred} shows the SBM predictions over a fixed horizon, scheduling applications rely on SBM predictions over a \textit{moving} horizon. 
Consider a scheduling horizon $t_\mathrm{sched} = 4$ days. 
At time $t=0$ days the SBM would be used to make predictions for days 0--4, and at time $t=1$, for days 1--5, etc.
In other words,  at any time $t$ we are interested in the instantaneous accuracy of the SBM predictions over the interval $[t, t+t_\mathrm{sched}]$, computed as:
\begin{equation}
\mathrm{MSE}_t = \frac{1}{n_\mathrm{sched}} \sum_{t'=t}^{t+n_\mathrm{sched}} (y_{t'} - \hat{y}_{t'})^2 \label{eq:mse}
\end{equation}
where $y_{t'}$ is the output value at time $t'$ \eqref{eq:arx}, $\hat{y}_{t'}$ is its predicted value, and $n_\mathrm{sched}$ is the number of samples taken during the interval $t_\mathrm{sched}$. 
Figure \ref{fig:sbm_mse} shows MSE$_t$ \eqref{eq:mse} for the same SBMs in Figure \ref{fig:sbm_pred} evaluated each hour for $t_\mathrm{sched} = 4$ days. The first seven days comprise the SBM training data, as in Figure \ref{fig:sbm_pred}. The vertical dashed lines mark times from $t=3$ days, where the four-day intervals first begin to include new data, to  $t=7$ days, where the four-day intervals comprise completely new data. 

For most of the SBMs, the four-day MSEs are relatively low for days 0--3. 
As data outside the training set are introduced (days 3--7), the four-day MSEs begin increasing.
The MSEs are highest past day seven, when all data are from outside the training set. 
Like the predictions in Figure \ref{fig:sbm_pred}, the moving-horizon prediction quality for these linear SBMs degrades after a few days, showing that re-fitting is needed for accurate predictions in scheduling calculations.

\section{Online Parameter Updating}

\subsection{Kalman Filter for Parameter Values}

\begin{figure}
	\begin{center}
		\includegraphics[width=\columnwidth]{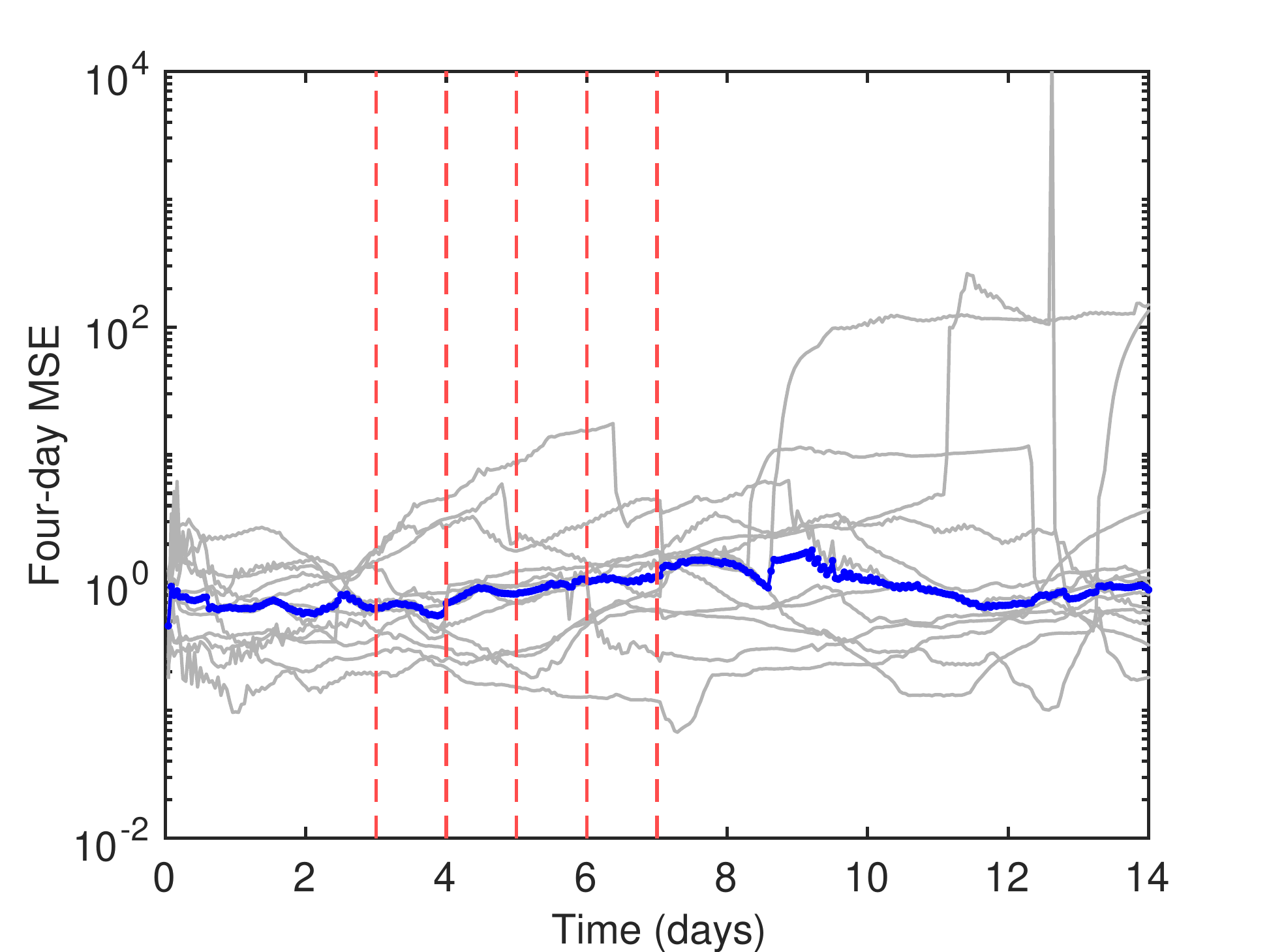}    
		\caption{SBM accuracy over a four-day moving horizon with online parameter estimation. Each gray line depicts one SBM, while the blue line indicates the median values. The dashed vertical lines mark days 3--7.}
		\label{fig:sbm_mse_kalman}
	\end{center}
\end{figure}

Along with least-squares methods, the Kalman filter is a well-established strategy for online parameter estimation \citep{ljung1990adaptation}, and its behavior and stability in this context have been analyzed by, e.g., \cite{cao2004analysis, guo1990estimating}.
Other works have also employed nonlinear extensions of the Kalman filter for online parameter estimation, e.g., the extended Kalman filter \citep{graichen2006feedforward}, or the unscented Kalman filter \citep{radecki2012online}. 

The underlying assumption is that the system state variables comprise the vector of model parameters, which vary via a random walk:
\begin{eqnarray}
\boldsymbol{\theta}_t = \boldsymbol{\theta}_{t-1} + \boldsymbol{\omega}_t \label{eq:randomwalk} \\
\boldsymbol{\omega}_t \sim N(\boldsymbol{0}, \boldsymbol{Q})
\end{eqnarray}
where $\boldsymbol{\theta}$ is the vector of parameters, and $\boldsymbol{Q}$ is the system noise covariance matrix. For the ARX model \eqref{eq:arx}, $\boldsymbol{\theta} = [-a_1,...,-a_N, b_{1,1},...,b_{1,N}, ..., b_{M,1},...,b_{M,N} ]^T$ and $\boldsymbol{\theta}_t$ is the estimate of $\boldsymbol{\theta}$ at time $t$. 
The model output is computed as a linear combination of the inputs: 
\begin{eqnarray}
y_t = \boldsymbol{u}_t^T \boldsymbol{\theta}_t + e_t \label{eq:arx2} \\
e_t \sim N(0, R) \label{eq:system2}
\end{eqnarray}
where $\boldsymbol{u}_t^T$ is the vector of inputs, and $R$ is the residual variance. For the SBM, $\boldsymbol{u}_t^T = [\phi_{1,t-1}, ... , \phi_{1,t-N}, \phi_{2,t-1},...,$ $\phi_{2,t-N}, \mathrm{SP}_{1,t-1}, ... , \mathrm{SP}_{1,t-N}, T_{t-1},...,T_{t-N}]$, making \eqref{eq:arx} and \eqref{eq:arx2} equivalent.
For this system, \cite{ljung1990adaptation} showed that a Kalman filter provides the best estimation of $\boldsymbol{\theta}_t$ when $\boldsymbol{Q}$ and $R$ are exactly known.  
Nevertheless, when $\boldsymbol{Q}$ and $R$ are unknown and/or the random walk \eqref{eq:randomwalk} does not model the true parameter dynamics well, the Kalman filter can often still track time-varying parameters accurately \citep{cao2004analysis, niedzwiecki2000identification}. 

The Kalman filter for \eqref{eq:randomwalk}--\eqref{eq:system2} can be written as:
\begin{eqnarray}
\hat{\boldsymbol{\theta}}_t = \hat{\boldsymbol{\theta}}_{t-1} + \boldsymbol{K}_t(y_t - \hat{y}_t) \label{eq:kalman1} \\
\hat{y}_t = \boldsymbol{u}_t^T \hat{\boldsymbol{\theta}}_{t-1} \\
\boldsymbol{K}_t = \boldsymbol{P}_{t-1} \boldsymbol{u}_k (R + \boldsymbol{u}_k^T \boldsymbol{P}_{t-1} \boldsymbol{u}_k)^{-1} \\
\boldsymbol{P}_t = [\boldsymbol{I} - \boldsymbol{K}_t \boldsymbol{u}^T] \boldsymbol{P}_{t-1} + \boldsymbol{Q} \label{eq:kalman2} 
\end{eqnarray}
where $\hat{\boldsymbol{\theta}}_t$ and $\hat{y}_t$ are the estimates for the parameter vector and the output, respectively, at time $t$. The vector $\boldsymbol{K}_t$ is the filter gain, and $\boldsymbol{P}_t$ is the estimated covariance matrix, both at time $t$. To maintain both nonzero $\boldsymbol{K}_t$ and stability, $\boldsymbol{P}_t$ should satisfy the standard matrix inequalities $\alpha \boldsymbol{I} \leq \boldsymbol{P}_t \leq \beta \boldsymbol{I}, \forall t$, where $\alpha$ and $\beta$ are positive scalars. 

We fix the residual variance $R=1$ for simplicity, leaving $\boldsymbol{Q}$ and $\boldsymbol{P}_0$ as tuning parameters for the filter. 
Assuming that the initial parameter estimates $\boldsymbol{\theta}_0$ are good at $t=0$, the initial covariance $\boldsymbol{P}_0$ should be small, and we set each parameter variance proportionally, such that $\boldsymbol{P}_0 = \mathrm{diag}( \boldsymbol{\theta}_0) \times 0.1 \%$.
For the system noise covariance $\boldsymbol{Q}$, which has less physical intuition, the parameter covariance matrix across the 12 training datasets, $\boldsymbol{\Sigma}$, was first computed. 
We then set $\boldsymbol{Q} = \boldsymbol{\Sigma} / n_\mathrm{train}$, where  $n_\mathrm{train}$ is the number of samples in the training dataset, such that the variance after $n_\mathrm{train}$ time steps is $\boldsymbol{\Sigma}$, or $(\boldsymbol{\theta}_t - \boldsymbol{\theta}_{t-n_\mathrm{train}}) \sim N(\boldsymbol{0}, \boldsymbol{\Sigma})$.

\subsection{Online Parameter Estimation Results}

The Kalman filter \eqref{eq:kalman1}--\eqref{eq:kalman2} was applied to the same SBM models using the one-minute process sampling interval. 
The initial parameter values $\boldsymbol{\theta}_0$ for each SBM were set to their respective estimated values from Section \ref{sec:sbm}. 
Figure \ref{fig:sbm_mse_kalman} shows the same hourly MSE$_t$ \eqref{eq:mse},  over a scheduling horizon $t_\mathrm{sched} = 4$ days, for the SBMs as in Figure \ref{fig:sbm_mse}, but with the SBM parameters updated at each one-minute sample via Kalman filter. 
For most of the SBMs, the four-day MSEs are again relatively low for days 0--3. 
While the four-day MSEs for some SBMs similarly begin increasing at day 3, the MSEs remain lower in general, owing to the updated parameter estimates.  
Unlike the SBMs with no Kalman filter in Figure \ref{fig:sbm_mse}, the median MSE remains relatively flat overall, showing that the online parameter estimation strategy effectively removes the need to periodically re-fit SBMs for use in scheduling calculations. 

\begin{figure}
	\begin{center}
		\includegraphics[width=0.88\columnwidth, trim = 0 1cm 0 0]{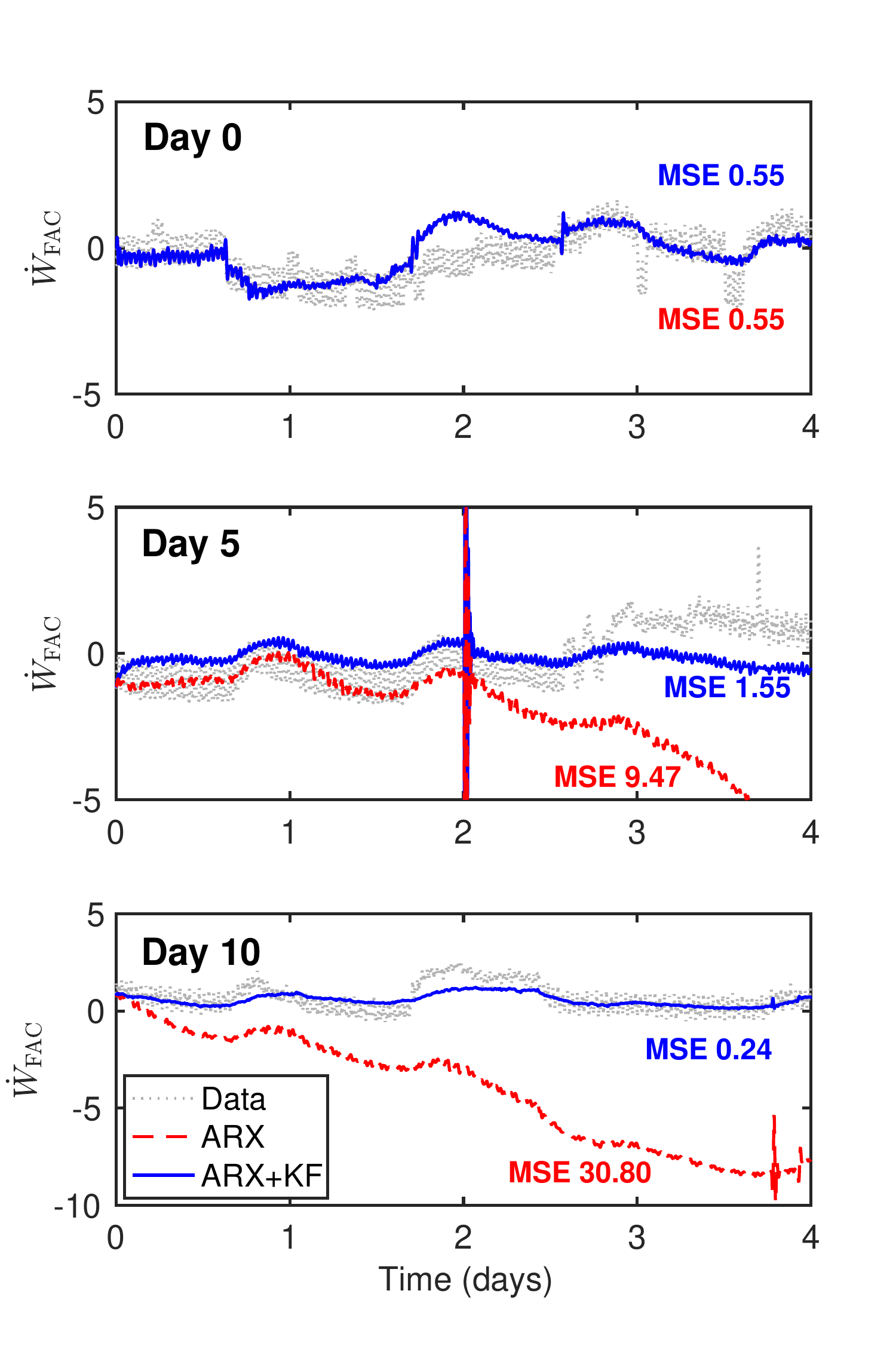}    
		\caption{Four-day predictions for ARX model with and without Kalman filter (KF) parameter updates at day zero (top), day five (middle), and day 10 (bottom).}
		\label{fig:compare}
	\end{center}
\end{figure}

The four-day predictions for the SBMs with and without parameter updates are shown in Figure \ref{fig:compare} for month five. 
These predictions correspond to the time points at zero, five, and ten days in Figures \ref{fig:sbm_mse} and \ref{fig:sbm_mse_kalman}.
Month five is shown because the SBM quickly loses accuracy without parameter updates (Figure \ref{fig:sbm_pred}). 
At day zero, the the Kalman filter has not updated the model parameters, and the two models are the same, with an MSE over the four-day window of 0.55. 
At day five, the time window comprises two days in the training dataset and two new days. 
Without parameter updates, the model is inaccurate for the two new days (Figure \ref{fig:compare}, middle), giving an overall MSE of 9.47.
With parameter updates, the model maintains an MSE of 1.55.
At day ten, when the four-day window does not include any training data, the updating strategy decreases the prediction MSE from 30.80 to 0.24. 
Note that the predictions at days five and ten differ from those in Figure \ref{fig:sbm_pred} because the initial conditions are updated when the models are evaluated in a moving horizon. 

\subsection{Parameter Dynamics}

\begin{figure}
	\begin{center}
		\includegraphics[width=\columnwidth]{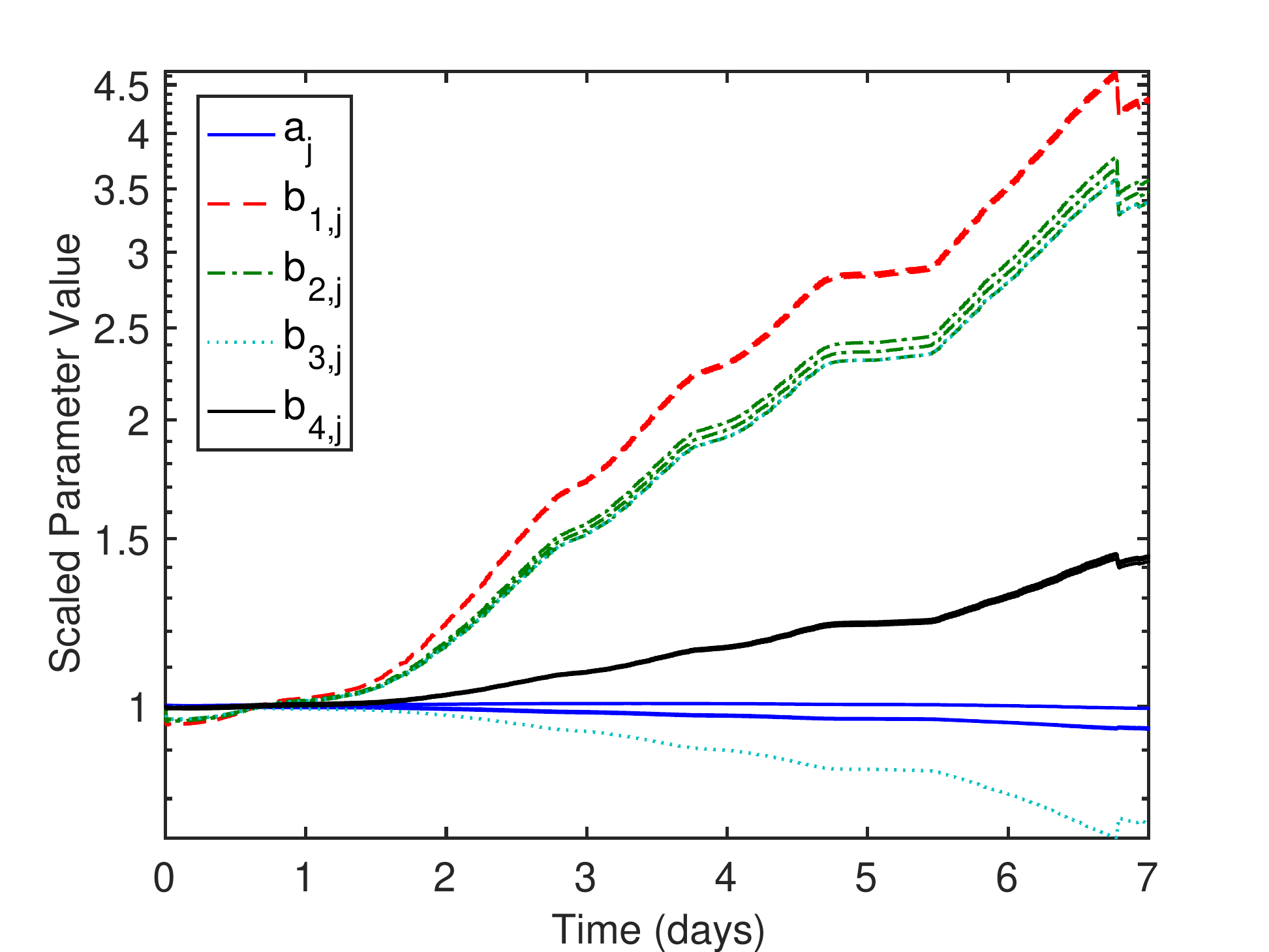}    
		\caption{Scaled SBM parameter estimates updated by Kalman filter over time. The three lines in each series depict three time-lagged coefficients.}
		\label{fig:paramupdates}
	\end{center}
\end{figure}

The values of the estimated parameter values may change significantly over time as they are updated by the Kalman filter \eqref{eq:kalman1}--\eqref{eq:kalman2}. 
These parameter ``dynamics'' are shown in Figure \ref{fig:paramupdates} for month five. 
As expected, the values of autoregressive parameters $a_j$ remain relatively constant in time. 
Since the historical data ($\boldsymbol{\Sigma}$) were used to tune the Kalman filter ($\boldsymbol{Q}$), the low (co)variances of $a_j$ in Figure \ref{fig:params_radar} result in low variations in their filtered dynamics. 
Similarly, the values of $b_{4,j}$ only increase by approximately a factor of 1.5, which is similar to their variation across monthly datasets (Figure \ref{fig:params_radar}). 
The values of $b_{1,j}$ and  $b_{2,j}$ increase relatively quickly, and future study should investigate the stability of SBM parameters with online estimation and/or alternative Kalman filter tuning(s). 
However, note that the four-fold increase is of a similar magnitude as the variations in $b_{1,j}$ and  $b_{2,j}$  shown in Figure \ref{fig:params_radar}. 
We verified that the parameters remained bounded for this process when the Kalman filter is applied over the entire year-long dataset.

While some parameters may approach new steady-state values (perhaps indicating a transition in the process dynamics), others may continue to fluctuate. 
For example, if the relationship between $y$ and $u_i$ is nonlinear, estimates of the parameters $b_{i,j}$ will be updated to approximate the local input-output response. 
Given the above, the dynamics of the SBM parameters could be monitored to reveal changes in the process/controller dynamics, or potentially aid in fault diagnoses. 
The dynamics of the SBM parameters may also reveal differences between ASU plants. 
For instance, in a transfer learning approach, the SBM parameters from one ASU could be used as the initial estimates $\boldsymbol{\theta}_0$ for another facility that has a similar design. 
Online parameter estimation can ``learn'' the SBM parameters of the new plant over time, revealing differences (or confirming similarities) between the plants.

\section{Conclusions}
Optimal demand-response scheduling of industrial ASUs should account for process dynamics and control to ensure that schedules are dynamically feasible. 
Here, scale-bridging models (SBMs) can represent closed-loop process dynamics in production scheduling using a low-order approximation. 
This work presents a data-driven strategy to create linear SBMs for an industrial process from routine operational data. 
We find that compact time-series models can accurately represent process dynamics over several days, but that their accuracy can suffer over time. 

To mitigate this issue, we apply Kalman filtering for online updating of the SBM parameters, including a tuning that incorporates the variations on the parameters in historical data. 
The strategy enables the SBM to adapt to nonlinear behavior and/or changes in the process dynamics, thereby greatly improving the accuracy of SBM predictions over four-day scheduling horizons. 
Furthermore, the dynamics of the parameter estimates can reveal changes in the process dynamics and potentially be used in monitoring applications. 
The Kalman filtering strategy can be used to adapt SBM to changes over extended periods of time, or to transfer an SBM from one plant to another. 

\section*{Disclaimer}
This report was prepared as an account of work sponsored by an agency of the United States Government. Neither the United States Government nor any agency thereof, nor any of their employees, makes any warranty, express or implied, or assumes any legal liability or responsibility for the accuracy, completeness, or usefulness of any information, apparatus, product, or process disclosed, or represents that its use would not infringe privately owned rights. Reference herein to any specific commercial product, process, or service by trade name, trademark, manufacturer, or otherwise does not necessarily constitute or imply its endorsement, recommendation, or favoring by the United States Government or any agency thereof. The views and opinions of authors expressed herein do not necessarily state or reflect those of the United States Government or any agency thereof.

\begin{ack}
We are grateful for support from CESMII Project \textit{Smart Manufacturing for Chemical Processing}, and the National Science Foundation  through CAREER  Award 1454433.
\end{ack}

\bibliography{ifacconf}             
\end{document}